\documentclass[11pt]{amsart}
\usepackage{amsmath}
\usepackage{amssymb}
\usepackage{amscd}
\usepackage{url}
\usepackage{pstricks}
\usepackage{pst-node}
\usepackage{pst-coil}
\usepackage{graphicx}

\newcommand{\R}{\mathbb{R}}
\newcommand{\intfracsmall}[2]{\genfrac{\lfloor}{\rfloor}{}{1}{#1}{#2}}
\newcommand{\intfraclarge}[2]{\left\lfloor\frac{#1}{#2}\right\rfloor}
\newcommand{\intfrac}[2]{\mathchoice{\intfraclarge{#1}{#2}}{\intfracsmall{#1}{#2}}
{\intfracsmall{#1}{#2}}{\intfracsmall{#1}{#2}}}
\newcommand{\rintfrac}[2]{\genfrac{\lceil}{\rceil}{}{1}{#1}{#2}}
\newcommand{\rintfraclarge}[2]{\left\lceil\frac{#1}{#2}\right\rceil}
\newcommand{\ceilfrac}[2]{\mathchoice{\rintfraclarge{#1}{#2}}{\rintfrac{#1}{#2}}
{\rintfrac{#1}{#2}}{\rintfrac{#1}{#2}}}
\newcommand{\intpart}[1]{\left\lfloor#1\right\rfloor}

\newcommand{\sawfracsmall}[2]{\genfrac{\langle}{\rangle}{}{1}{#1}{#2}}
\newcommand{\sawtooth}[1]{\left\langle #1\right\rangle}
\newcommand{\sawfraclarge}[2]{\sawtooth{\frac{#1}{#2}}}
\newcommand{\sawfrac}[2]{\mathchoice{\sawfraclarge{#1}{#2}}
{\sawfracsmall{#1}{#2}}{\sawfracsmall{#1}{#2}}{\sawfracsmall{#1}{#2}}}

\newcommand{\bp}{\begin{pmatrix}}
\newcommand{\ep}{\end{pmatrix}}

\newcommand{\epn}{\varepsilon} 

\newcommand{\calt}{{\mathcal T}}

\numberwithin{equation}{section}

\theoremstyle{plain}
\newtheorem{theorem}[equation]{Theorem}
\newtheorem{lemma}[equation]{Lemma}
\newtheorem{proposition}[equation]{Proposition}

\newtheorem{corollary}[equation]{Corollary}

\theoremstyle{definition}

\newtheorem{remark}[equation]{Remark}

\newtheorem*{acknowledgements}{Acknowledgements}

\def\Z{\mathbb Z}\def\Q{\mathbb Q}\def\N{\mathbb N}
\def\S{{\mathcal S}_{pq}}
\def\pone{$(+1)$}
\def\mone{$(-1)$}

\title[surgeries on torus knots]{Heegaard Floer homologies of \pone\ surgeries on torus knots}
\author{Maciej Borodzik}
\address{Institute of Mathematics, University of Warsaw, ul. Banacha 2,
02-097 Warsaw, Poland}
\email{mcboro@mimuw.edu.pl}

\author{Andr\'as N\'emethi}
\address{A. R\'enyi Institute of Mathematics, 1053 Budapest,  Re\'altanoda u. 13-15,  Hungary.}
\email{nemethi@renyi.hu}

\date{\today}
\makeatletter
\@namedef{subjclassname@2010}{\textup{2010} Mathematics Subject Classification}
\makeatother
\begin{document}
\begin{abstract}
We compute the Heegaard Floer homology of $S^3_1(K)$
(the \pone\ surgery on the torus knot $T_{p,q}$)
in terms of the semigroup generated by $p$ and $q$, and
we find a compact formula (involving Dedekind sums) for the corresponding Ozsv\'ath--Szab\'o $d$-invariant.
We relate the result  to
known knot invariants of $T_{p,q}$ as the genus and the Levine--Tristram signatures. Furthermore,
 we emphasize the striking resemblance between Heegaard Floer homologies
of \pone\ and \mone\ 
surgeries on torus knots. This relation is best seen at the level of  $\tau$ functions.
\end{abstract}
\maketitle

\section{Introduction}

Let $K\subset S^3$ be a knot. Let us consider $S^3_1(K)$, the \pone\ surgery of $K$.
The main goal of the present article is the determination of the
Heegaard Floer homology $HF^+$ with $\Z$--coefficient of $S^3_1(K)$, when
$K=T_{p,q}$ is the torus knot.
Note that $S^3_1(K)$ is an integral homology sphere, hence
its Heegaard Floer homology is concentrated in 
its unique $spin^c$--structure.
Moreover, we provide several closed formulae
for the correction term
$d(K):=d(S^3_1(K))$  of Ozsv\'ath--Szab\'o \cite{OS}. In fact, searching for such closed formulae,
and the recent article of Peters \cite{Pet} regarding several properties of $d(K)$ (see Section~\ref{sec:D-prop} below)
was the motivation and the starting point of the present work.

\subsection{The Heegaard Floer homology.} Let us fix two relative prime integers $p$ and $q$. We set
$\delta:=(p-1)(q-1)/2$ and we define $\S$ as the subsemigroup of $\N$ generated by $p$ and $q$ including 0 too.
It is well--known that $\N\setminus \S$ is finite of cardinality $\delta$, hence the integers
\begin{equation}\label{eqint:1}
\alpha_i:= \#\{s\not\in\S\,:\, s>i\}
\end{equation}
are well defined. In fact, $\delta=\alpha_0\geq \alpha_1\geq \cdots \geq \alpha_{2\delta-2}=1$,
and $\alpha_i=0$ for $i>2\delta-2$. These integers --- or, equivalently, the semigroup $\S$ --- codifies
the same amount of data as the Alexander polynomial $\Delta(t)$ of the knot, or the symmetrized Alexander
polynomial $\Delta^{\#}(t)=t^{-\delta}\Delta(t)$ (see Section \ref{section:hf}).

In the description of the Heegaard  Floer homology we use (the already standard) notations of $\Z[U]$--modules,
what is also recalled in Section \ref{section:hf}.
\begin{theorem}\label{thmintr:HF} For the Seifert 3--manifold $-\Sigma=S^3_1(T_{p,q})$ one has
\begin{align*}
&HF^+_{odd}(-\Sigma)=0,\ \ \
HF^+_{even}(-\Sigma)=\calt^+_{d(-\Sigma)}\oplus \bigoplus_{k=0}^{\delta-2}
\calt_{k(k+1)-2\alpha_{\delta+k}}(\alpha_{\delta+k})^{\oplus  2},\\
&d(-\Sigma)=-2\alpha_{\delta-1},\ \ \ \mbox{and} \ \ \
{\rm rank}_{\Z}\,
HF^+_{red}(-\Sigma)=\frac{1}{2}(\Delta^{\#})''(1)-\alpha_{\delta-1}.
\end{align*}
\end{theorem}

\begin{remark}\label{rem:intr} (a) If we write $\Delta^{\#}(t)=a_0+\sum _{j=1}^{\delta}a_j(t^j+t^{-j})$, then using the
correspondence between the Alexander polynomial and the semigroup recalled  in  Section \ref{section:hf}, one can
show that $\alpha_{\delta-1}=\sum_{j\geq 1}ja_j$. In this way we recover the formula of
Ozsv\'ath and Szab\'o \cite[Section~7]{OS-al}: $d(S^3_1(K))=-2\sum _{j\geq 1}ja_j$, but this time in terms of the
semigroup $\S$, a fact which provides a new geometric interpretation. Note that the coefficients $\alpha_i$
enter in a substantial way in the reduced  part of $HF^+$ too.

(b) Since the Casson invariant $\lambda $ (or, the Seiberg--Witten invariant ${\mathbf {sw}}$)
of $-\Sigma$ satisfies $\lambda={\mathbf{ sw}}= {\rm rank}HF^+_{red}-d/2$, we get
$\lambda(-\Sigma)={\mathbf{sw}}(-\Sigma)=\frac12(\Delta^{\#})''(1)$, the classical surgery formula
for $\lambda$, see \cite{Le} and \cite{AMC}.
\end{remark}

\subsection{The $d$--invariant.}\label{sec:D-prop} For any 3--manifold $M$, the collection of correction terms $d(M,\sigma)$
associated with its $spin^c$--structures capture very strong geometric information, and recently they
provided many deep applications. If $K$ is a knot,  and $M=S^3_1(K)$, then
$d(K):=d(S^3_1(K))$
was intensively studied  by Peters \cite{Pet}, who showed, among other things, that $d(K)$ is a
concordance invariant and it provides  a bound for four genus:
\begin{equation}\label{eq:fourgenus}
0\le -d(K)\le 2g_4(K).
\end{equation}
Furthermore, if $K'$ arises from $K$ by changing one \emph{negative} crossing into a
\emph{positive} one on some diagram, one has
the following relation
\begin{equation}\label{eq:unknot}
d(K)-2\le d(K')\le d(K).
\end{equation}
It follows, in particular, that $-\frac12(d(K)+d(mK))\le u(K)$, where
$u(K)$ is the unknotting number and $mK$ is the mirror of $K$.
Thus $d(K)$ bears  a strong resemblance to the classical signature $\sigma(K)$ of $K$.
Indeed, by \cite{OS-al}, if $K$ is alternating
then $d(K)=2\min\left(0,\lfloor\frac{\sigma(K)}{4}\rfloor\right)$.

\vspace{1mm}

In the present article we provide two further
 closed formulae of $d(T_{p,q})$. The first one is in terms of generalized Dedekind sums,
 the second one is in terms of Levine--Tristram signature of the torus knot  evaluated at
$\exp(2\pi i\frac{\delta}{pq})$.
\begin{theorem}\label{thm:main}
Set $c:=0$ if $\delta-1\in\S$, otherwise take $c:=1$.
Furthermore, for two integers $a$ and $b$, $a\neq 0$ we write
\[
\epn_{a}(b)=
\begin{cases}
1&\mbox{if} \ \ a|b\\
0&\text{otherwise, }
\end{cases} \ \ \ \mbox{and} \ \  \ \epn_{p,q}(b)=\epn_p(b)+\epn_q(b).
\]

Then
\begin{multline*}
d(S^3_1(T_{p,q}))=-1-\frac{(\delta-1)(\delta-2)}{pq}-\frac{1}{6pq}(p^2+q^2-3p-3q-2)-\\-\ceilfrac{\delta}{p}
-\ceilfrac{\delta}{q}
+c+\frac12\epn_{p,q}(\delta-1)+2s(p,q\,;\frac{\delta-1}{q},0)+2s(q,p\,;\frac{\delta-1}{p},0).
\end{multline*}
where $s(p,q\,;x,0)$ are the generalized Dedekind sums (see Section~\ref{section:preliminary} below).

(For the second  formula see Section   \ref{section:TL}.)
\end{theorem}

For $q=p+1$ we obtain 
$$d(S^3_1(T_{p,p+1}))=-\intpart{\frac{p}{2}}\big(\intpart{\frac{p}{2}}+1\big).$$
For $p=2$ we get $d(S^3_1(T_{2,2\delta+1}))=-2\ceilfrac{\delta}{2}$; since in this case $\sigma(K)=-2\delta$,
this identity is compatible with the above formula of Ozsv\'ath and Szab\'o valid for alternating knots.

Our approach permits us to find sharp inequalities for $d(K)$ valid for $K=T_{p,q}$.  Recall that in this case
$g_4=\delta$ \cite{KM}.
\begin{corollary}\label{cor:INEQ} (a)
$$-d(S^3_1(T_{p,q}))\leq g_4+1\leq -\sigma(K).$$
The first inequality is equality for $p=2$ and $\delta $ odd.

(b) Assume that $p<q$ and $2q-\frac52\leq c\sqrt{\delta}$ for some constant $c>0$. Then
$$-d(S^3_1(T_{p,q}))\leq 2q-2+\frac{\delta-1}{2}
\leq \frac{\delta}{2}+c\sqrt{\delta}\leq \frac{-\sigma-1}{2}+c\sqrt{-\sigma}.$$
\end{corollary}
Part (b) is `strong' whenever  $p$ `grows together with $q$'. In this case we reobtain at least  asymptotically
the growth $-\sigma/2$ valid for alternating knots.

We expect similar inequalities  for $d(S^3_1(K))$ for any algebraic knot $K$.

\subsection{The methods and the structure of the article.}
The main tool of the proof is the $\tau$ function associated with any Seifert (or `almost rational')
negative definite plumbed manifold  \cite{Nem}. For this,
 we interpret $S^3_1(T_{p,q})$ as a Seifert manifold and
run the algorithm of [loc.cit.].
We study the corresponding $\tau$ function in the same way, as in \cite[Section~7]{Nem-neg2}
and are able to obtain its local maxima and minima (cf. Section~\ref{section:proof}).
This data provides the corresponding `graded roots' (cf. \cite{Nem,Nem-grad}), hence all the information
regarding the Heegaard Floer homology (see Section~\ref{section:hf}).

We recall that Peters studies $d_-(K)=d(S^3_{-1}(K))$ too. But, for any positive knot, $d_-(K)=0$.
(See  also  Proposition \ref{prop:s'}(b) below.)
[This also follows  from the fact $d_-(K)=d(mK)$. Indeed,
as $mK$ can be unknotted by changing only negative crossings to positive ones,  by \eqref{eq:unknot}
 we get $d(mK)\ge d(U)=0$, where $U$ is
the unknot.  But, by \eqref{eq:fourgenus}, $d(mK)\le 0$ as well.]

Although the qualitative behavior of $S^3_{-1}(T_{p,q})$ and $-S^3_1(T_{p,q})$ at the level of $d$--invariants
is different,  we wish  to stress that
the $\tau$ functions related to these two manifolds
turn out to be strongly related to each other. This relation can be clearly seen at the level of $\tau$ functions,
 but becomes
less transparent when we pass to graded roots or the Heegaard Floer homologies. For more details, see
Section \ref{section:hf}.
We are wondering, whether this `duality'
holds for a larger group of knots, e.g. for all algebraic knots.
[As for general algebraic knot $K$,  neither $S^3_1(K)$
nor $-S^3_1(K)$ has negative definite plumbing graph, a possible generalization of
 this `duality'  might use an  extension of  the results of  \cite{Nem-neg2}
for not necessarily negative definite graphs.]

In Section~\ref{section:TL} we relate the invariants to Levine--Tristram signatures.
Finally, in the last section we establish the  inequalities of Corollary \ref{cor:INEQ}.

\section{Preliminary results}\label{section:preliminary}
First we recall some results from \cite[Section 11]{Nem}. Let $\Sigma=\Sigma(e_0,(\alpha_1,\omega_1),\dots,(\alpha_\nu,\omega_\nu))$
be a plumbed negative definite three--manifold. The notation means that the plumbing diagram
is star-shaped, the central vertex has weight $e_0$
and there are $\nu$ arms stemming from it.
Furthermore, $0<\omega_l<\alpha_l$, $\gcd(\alpha_l,\omega_l)=1$ ($1\leq l\leq \nu$),  and
if we write  the continued  fraction
\[\frac{\alpha_l}{\omega_l}=k_{l1}-\cfrac{1}{k_{l2}-\cfrac{1}{\ddots-\cfrac{1}{k_{ls_l}}}}
\] with all $k_{ln}\geq 2$
then the $l$--th arm consists 
of a chain of $s_l$ vertices with weights $-k_{l1}$,\dots,$-k_{ls_l}$. The weight of the vertex closest to the central one is $-k_{l1}$.

Let us introduce the quantities
\[
e=e_0+\sum_{l=1}^\nu\frac{\omega_l}{\alpha_l}
\text{\ \ and \ \ }
\varepsilon=\left(2-\nu+\sum_{l=1}^\nu\frac{1}{\alpha_l}\right)\frac{1}{e},
\]
and the notation $\alpha:=\prod_l\alpha_l$ and $\hat{\alpha}_l:=\alpha/\alpha_l$.
We assume that $e<0$. Observe that
\[|H_1(\Sigma,\mathbb{Z})|=-e\alpha_1\dots\alpha_\nu,\]
in particular, $\Sigma$ is an integral homology sphere if and only if $e=-\frac{1}{\alpha_1\dots\alpha_\nu}$, that is
\begin{equation}\label{eq:inthom}
-1=e_0\alpha+\sum_{l}\omega_l\hat{\alpha}_l.
\end{equation}
Hence the integers $\omega_l$'s are determined by the $\alpha_l$'s.
Let us consider the canonical  $spin^c$--structure $\sigma$ on $\Sigma$. We have the following fact.
\begin{proposition}[see \expandafter{\cite[Theorems 11.9, 11.12 and 11.13]{Nem}}]
The $d$ invariant of $\Sigma=(\Sigma,\sigma)$ is given by the following formula
\begin{equation}\label{eq:main}
d(\Sigma)=\frac{1}{4}(K^2+s)-2\min_{m\ge 0}\tau(m),
\end{equation}
where
\begin{equation}\label{eq:A}
K^2+s=\varepsilon^2e+e+5-12\sum_{l=1}^\nu s(\omega_l,\alpha_l),
\end{equation}
and $\tau(m)$ is a function defined by
\begin{equation}\label{eq:B}
\tau(m)=\sum_{j=0}^{m-1}\Delta_j, \ \ \mbox{where} \ \ \
\Delta_j:=1-je_0-\sum_{l=1}^\nu\ceilfrac{j\omega_l}{\alpha_l}.
\end{equation}
\end{proposition}
\noindent $K^2+s$ is an invariant of a plumbed 3--manifold computed from its
plumbing graph: $K^2$ is the self--intersection of the
canonical class and $s$ is the number of
vertices of the graph; however, in this paper we only use formula~\eqref{eq:A}.
Note that the Heegaard Floer homology is also determined in terms of the $\tau$--function.
For details, see e.g. \cite{Nem,Nem-neg}.

\vspace{1mm}

We recall also a definition of the sawtooth function: for  $x\in\mathbb{R}$ one sets
\[\sawtooth{x}=\begin{cases}
0&x\in\Z\\
x-\intpart{x}-\frac12&x\not\in\Z.
\end{cases}\]
(We use the notation $\sawtooth{x}$ and not $((x))$, because we think it is then easier to read formulae).
The generalized Dedekind sum is defined for two integers, $p,q\in\Z\setminus\{0\}$
and $x,y\in\mathbb{R}$ by
\begin{equation}\label{eq:s-def}
s(p,q;x,y)=\sum_{i=0}^{q-1}\sawtooth{\frac{i+y}{q}}\sawtooth{p\frac{i+y}{q}+x}.
\end{equation}
The classical Dedekind sum is
$s(p,q)=s(p,q;0,0)$, cf.  \cite{RG}.
We introduce one  more notation. Recall that $\epn_a(b)$ was defined in the statement of
 Theorem \ref{thm:main}. For a set of Seifert invariants  $\{\alpha_l\}_l$ we also write
$\epn_{\{\alpha\}}(b):=\sum_l\epn_{\alpha_l}(b)$.

\section{A compact formula for $\tau(i)$}\label{section:compact}
Under conditions $H_1(\Sigma,\Z)=0$   we provide a compact formula for $\tau(m)$.
\begin{proposition}\label{prop:tauform}
Let us write $m=d(m)\alpha+\sum_la_l\hat{\alpha}_l$, for some integers $d(m)$ and  $0 \leq a_l <\alpha_l$ for all $l$. Then
\begin{equation}\label{eq:tau-formula}
\begin{split}
\tau(m)=&\sum_l\left(\frac12 \intfrac{m-1}{\alpha_l}
-s(\hat{\alpha}_l,\alpha_l\,;\frac{m}{\alpha_l},0)+s(\hat{\alpha}_l,\alpha_l)\right)\\
&+\frac{m^2}{2\alpha}+m(1-\frac{\nu}{2})-\frac{d(m)}{2}+\frac{\nu}{4}+\frac{1}{4}\epn_{\{\alpha\}}(m).
\end{split}
\end{equation}
(For a slightly different formula --- using the $\omega_l$'s but no $d(m)$ --- see (\ref{eq:difver}) in the proof.)
\end{proposition}

\begin{proof}
After straightforward computations using \eqref{eq:inthom} we get
\begin{equation}\label{eq:Deltak}
\Delta_j=\sum_l\sawfrac{j\omega_l}{\alpha_l}+1-\frac{\nu}{2}+\frac{j}{\alpha}+\frac{\epn_{\{\alpha\}}(j)}{2}.
\end{equation}
We need a following result whose proof is at the end of this section.
\begin{lemma}\label{lem:sum}
Let $a$ and $b$ be two coprime integers, with $b\neq 0$. Then, for any $m\in\mathbb{Z}$,
\begin{equation}\label{eq:lemsum}
\sum_{j=0}^{m-1}\sawfrac{ja}{b}=-\frac12\sawfrac{ma}{b}-s(a^*,b\,;\frac{m}{b},0)+s(a^*,b),
\end{equation}
where $a^*$ is defined by the condition  $a^*a\equiv-1 \bmod b$.
\end{lemma}
\noindent By this lemma applied for each  $\sawfrac{j\omega_l}{\alpha_l}$ in \eqref{eq:Deltak}, and using
  $\omega_l^*=\hat{\alpha}_l\bmod \alpha_l$, cf. \eqref{eq:inthom},
 \begin{equation}\label{eq:difver}
\begin{split}
\tau(m)=&\sum_l\left(-\frac12\sawfrac{m\omega_l}{\alpha_l}+\frac12 \intfrac{m-1}{\alpha_l}
-s(\hat{\alpha}_l,\alpha_l\,;\frac{m}{\alpha_l},0)+s(\hat{\alpha}_l,\alpha_l)\right)\\
&+m(1-\frac{\nu}{2})+\frac{\nu}{2}+\frac{m(m-1)}{2\alpha}.
\end{split}\end{equation}
Then using the sum of
$$\sawfrac{m\omega_l}{\alpha_l}=\sawfrac{a_l\hat{\alpha}_l\omega_l}{\alpha_l}
=\sawfrac{-a_l}{\alpha_l}=-\frac{a_l}{\alpha_l}+\frac{1}{2}-\frac{1}{2}\epn_{\alpha_l}(m)$$
we conclude the proof.
\end{proof}
\begin{proof}[Proof of Lemma~\ref{lem:sum}]
Since
\begin{equation}\label{eq:dag}
\sum_{j=0}^{b-1}\sawfrac{ja}{b}=0
\end{equation} (see e.g. \cite[Lemma 1]{RG})
we may assume that $0\le m<b$. Then we can write
\[\sum_{j=0}^{m-1}\sawtooth{\frac{ja}{b}}=\sum_{j=0}^{m-1}\sawtooth{\frac{ja}{b}}
\ceilfrac{m-j}{b}=\sum_{j=0}^{b-1}\sawtooth{\frac{ja}{b}}
\ceilfrac{m-j}{b}.\]
Upon using the definition of $\sawtooth{\cdot} $  the last sum becomes
\[-\frac12\sawtooth{\frac{ma}{b}}+\sum_{j=0}^{b-1}\sawtooth{\frac{ja}{b}}
\left(\frac12+\sawtooth{\frac{-m+j}{b}}+\frac{m-j}{b}\right).\]
Using \eqref{eq:dag} again,  we can rewrite this as
\[
-\frac12\sawfrac{ma}{b}+\sum_{j=0}^{b-1}
\sawfrac{ja}{b}\sawfrac{-m+j}{b}-\sawfrac{ja}{b}\sawfrac{j}{b}.
\]
Now, the substitution $j\mapsto -a^*j$ and identity $\sawtooth{-x}=-\sawtooth{x}$ provide the statement.
\end{proof}

\section{Extrema of the function $\tau(i)$}\label{section:proof}

It is well-known (see e.g. \cite{Mos}), that if $K=T_{p,q}$ is a torus knot (for $p$ and $q$ relative prime
positive integers), then
\[S_{+1}(K)=-\Sigma(e_0,(p,p'),(q,q'),(r,r')),\]
where $r=pq-1$, $p'$, $q'$, $r'$ (the $\omega_l$'s) and $e_0$ satisfy \eqref{eq:inthom}.
In fact, $e_0=-2$ and  $p'$, $ q'$ and $r'$ are determine uniquely by
\[p'q\equiv 1\bmod p,\ pq'\equiv 1\bmod q,\ r'=pq-2.\]

The minus sign in front of $\Sigma$ shows
the change of orientations.
We also write $\delta:=(p-1)(q-1)/2$.
Using \eqref{eq:A} via  $s(p',p)=s(q,p)$, the Dedekind reciprocity law and  a direct computation we get for $\Sigma$
\begin{equation}\label{eq:K}(K^2+s)(\Sigma)=-4\delta(\delta-3).\end{equation}
In this section we study the local extrema of the function $\tau(i)$
associated with  $\Sigma=\Sigma(e_0,(p,p'),(q,q'),(r,r'))$. We find also the global minimum, which together with
results from Section~\ref{section:compact} gives the proof of Theorem~\ref{thm:main}.

\begin{proposition}\label{prop:minoftau}

Let $\S$ be the subsemigroup of ${\mathbb N}$  generated by $p$ and $q$ and including 0. Then the following facts hold:

(a) \
The function $\tau:{\mathbb N}\to \Z$ attains its local minima at values $m_n=n(pq-1)$ for $0\leq n\leq 2\delta-2$,
and local maxima at values $M_n=npq+1$ for $0\leq n\leq 2\delta-3$. This means that $\tau$ is (not necessarily
strict) increasing on any interval $[m_n,M_n]$ and $[m_{2\delta-2},\infty)$, and  (not necessarily
strict) decreasing on any interval $[M_n,m_{n+1}]$.

(b) \ The sequences $\{m_n\}_n$ and $\{M_n\}_n$ are minimal with these properties, that is:
$$\tau(M_n)>\max\{\tau(m_n),\tau(m_{n+1})\}.$$ In fact, for any $0\leq n\leq 2\delta-3$, one has
\begin{equation}\label{eq:SEMI}
\begin{split}
\tau(M_n)-\tau(m_{n+1})=& \ \#\{s\not\in\S\,:\, s\geq n+2\}>0,\\
\tau(M_{n})-\tau(m_{n})=& \ \#\{s\in\S\,:\, s\leq n\}>0.
\end{split}
\end{equation}

(c) \ The absolute minimum of $\tau$ occurs for $m_{\min}=m_{\delta-1}$.

(d) \ For any $0\leq n\leq 2\delta-3$,
\begin{equation}\label{eq:Mn}
\tau(M_n)=\frac{n(n-2\delta+3)}{2}+1.\end{equation}
\end{proposition}
\begin{proof} First note that $\tau(m_0)=\tau(0)=0$ while $\tau(M_0)=\tau(1)=\Delta_0=1$.

We define  $M_{n}=npq+1$ for any $n$ and
we will compute $\tau(j+1)-\tau(j)=\Delta_j$ for any $j\in [M_n,M_{n+1})$. Clearly, $M_n\leq j < M_{n+1}$
 if and only if  $0\leq (n+1)pq-j<pq$.

Note the following fact regarding the semigroup $\S$ and any integer $a\in[0,pq)$:
\begin{equation}\label{eq:sgr}
\begin{split}
a\in \S \ & \Longleftrightarrow \ \ \ \ \ \ \  a=\alpha p+\beta q \ \ \ \ (0\leq \alpha<q,\ 0\leq \beta< p),\\
a\not\in \S \ & \Longleftrightarrow \ a+pq=\alpha p+\beta q \ \ \ \ (0\leq \alpha<q,\ 0\leq \beta< p).
\end{split}
\end{equation}
First we fix an interval $[M_n,M_{n+1})$ for some $0\leq n\leq 2\delta-3$.

\vspace{1mm}

\noindent {\bf Case 1.} Assume that $(n+1)pq-j\in\S$. Then, by \eqref{eq:sgr}
$$j=(n+1)pq-\alpha p-\beta q \ \ \ \ (0\leq \alpha<q,\ 0\leq \beta< p).$$
In particular, $jp'\equiv -\beta\bmod p$, $jq'\equiv -\alpha\bmod q$, and in general, $jr'\equiv -j\bmod r$.
Using $-\sawtooth{-\beta/p}=\sawtooth{\beta/p}=\beta/p-1/2+\epn_p(\beta)/2=\beta/p-1/2+\epn_p(j)/2$,
and similarly for $\alpha/q$, by a computation the value of $\Delta_j$ from  \eqref{eq:Deltak}
transforms into $\Delta_j=\intpart{j/r}-n$. Therefore
\begin{equation}\label{eq:Delta1}
\Delta_j=\intpart{\frac{j}{r}}-n=\left\{
\begin{array}{ll}0 & \mbox{if $M_n\leq j<m_{n+1}$}\\
1 & \mbox{if $m_{n+1}\leq j< M_{n+1}$}.\end{array}\right.\end{equation}

\noindent {\bf Case 2.} Assume that $(n+1)pq-j\not\in\S$. Then, by \eqref{eq:sgr}
$$j=(n+2)pq-\alpha p-\beta q \ \ \ \ (0\leq \alpha<q,\ 0\leq \beta< p).$$
Similar computation as in in case 1 provides
\begin{equation}\label{eq:Delta2}
\Delta_j=\intpart{\frac{j}{r}}-n-1=\left\{
\begin{array}{ll}-1 & \mbox{if $M_n\leq j<m_{n+1}$}\\
0 & \mbox{if $m_{n+1}\leq j< M_{n+1}$}.\end{array}\right.\end{equation}

\noindent {\bf Case 3.} Assume that $j\in [M_n,M_{n+1})$ with $n\geq 2\delta-2$. If $(n+1)pq-j\in\S$,
then by a similar argument as in Case 1 we get  $\Delta_j=\intpart{j/n}-n$. But $j\geq M_n\geq nr$, hence
$\Delta_j\geq 0$.

If $(n+1)pq-j\not\in\S$, then $(n+1)pq-j\leq 2\delta-1$ (the conductor of $\S$), hence automatically
$j\geq (n+1)pq-2\delta+1$. But $(n+1)pq-2\delta+1\geq (n+1)r$ whenever $n\geq 2\delta-2$. Hence
$\Delta_j=\intpart{j/r}-n-1\geq 0$ again. This ends the proof of (a).

Moreover, equations \eqref{eq:Delta1} and \eqref{eq:Delta2} provide \eqref{eq:SEMI} too.
Since $2\delta-1\not\in\S$ and $0\in\S$, both differences are strictly positive. This proves (b).

\smallskip
In order to prove (c), let us find the difference between two subsequent local minima.
Note that $s\in\S$ if and only if $2\delta-1-s\not\in\S$. Therefore, via \eqref{eq:SEMI},
the difference  equals
$$\tau(m_{n+1})-\tau(m_n)=\ \#\{s\not\in\S\,:\, s\geq 2\delta-1-n\}-\#\{s\not\in\S\,;\, s\geq n+2\}.$$

\noindent For $n\leq \delta-2$
one gets $2\delta-1-n\geq n+2$, hence this difference is $\leq 0$. Similarly,
$$\tau(m_{n+1})-\tau(m_n)=\ \#\{s\in\S\,:\, s\leq n\}-\#\{s\in\S\,;\, s\leq 2\delta-1-(n+2)\},$$
which is $\geq 0$  for $n\geq \delta-1$. This ends part (c) too.

Part (d) can be verified in two ways. First, by \eqref{eq:SEMI} we deduce that $\tau(M_{m+1})-\tau(M_n)$
is
$$ \#\{s\in\S\,:\, s\leq n+1\}-\#\{s\not\in \S\,:\, s\geq n+2\}=-\#\{s\not\in\S\}+\#\{0\leq s\leq n+1\}$$
which is $n+2-\delta$. Then one proceeds by induction. Alternatively, by a direct computation one gets that
$\tau(npq)=n(n-2\delta+3)/2 $ (here the Dedekind reciprocity law is used) and $\Delta_{npq}=1$.
\end{proof}

\section{Proof of Theorem~\ref{thm:main}}\label{section:prooftrue}
Having studied the monotonicity properties of the $\tau$-function, we compute the value $\tau(m_{\min})$ via \eqref{eq:difver}.
By Proposition~\ref{prop:minoftau}, $m_{\min}=(\delta-1)r$.
We write  $\delta-1=ap+bq-cpq$, with $0\le a<q$, $0\le b<p$ and $c\in\{0,1\}$; cf.
\eqref{eq:sgr}.  We have
\begin{equation*}
\intfrac{m-1}{p}+\intfrac{m-1}{q}+\intfrac{m-1}{r}
=(\delta-1)(p+q)-\ceilfrac{\delta}{p}-\ceilfrac{\delta}{q}+\delta-2.
\end{equation*}
Moreover
$\sawfrac{p'm}{p}=\sawfrac{-p'(\delta-1)}{p}=\sawfrac{-bp'q}{p}=\sawfrac{-b}{p}$, and $\sawfrac{q'm}{q}=\sawfrac{-a}{q}$.
Their sum is
\[\sawfrac{-a}{q}+\sawfrac{-b}{p}=1-\frac{a}{q}-\frac{b}{p}-\frac12\epn_{p,q}
(\delta-1)=1-\frac12\epn_{p,q}(\delta-1)-\frac{\delta-1}{pq}-c.\]

As for the Dedekind sums, we observe that $s(pq,r\,;\frac{m}{r},0)=s(pq,r)$ so these terms cancel each other in \eqref{eq:difver}.  
Moreover, $s(pr,q\,;\frac{m}{q},0)=s(-p,q\,;-\frac{\delta-1}{q},0)=-s(p,q\,;\frac{\delta-1}{q},0)$, and
by the reciprocity law
\[s(p,q)+s(q,p)=\frac{1}{12pq}(p^2+q^2+1-3pq)=\frac{1}{12pq}(p^2+q^2-3p-3q-2)-\frac{\delta-1}{2pq}.\]
Putting this together we get
\begin{multline*}
2\tau(m_{\min})=-\delta^2+3\delta-1+\frac{-\delta^2+3\delta-2}{pq}-\frac{1}{6pq}(p^2+q^2-3p-3q-2)\\
-\ceilfrac{\delta}{p}-\ceilfrac{\delta}{q}
+c+\frac12\epn_{p,q}(\delta-1)+2s(p,q\,;\frac{\delta-1}{q},0)+2s(q,p\,;\frac{\delta-1}{p},0).
\end{multline*}
Hence, using \eqref{eq:main}, \eqref{eq:K}
 and $d(S^3_1(T_{pq}))=-d(\Sigma)$ we obtain  formula (a) of  Theorem~\ref{thm:main}.
Formula (b) of \ref{thm:main} follows from the first equation of \eqref{eq:SEMI} written for
$n=\delta-2$, \eqref{eq:Mn} and \eqref{eq:K}.

\section{The Heegaard Floer homology of $\Sigma=-S^3_{1}(T_{p,q})$}\label{section:hf}

Let $\Delta(t)$ be the Alexander polynomial of the knot $K=T_{p,q}$ normalized in such a way that
$\Delta(1)=1$. One has the following identity  connecting
$\Delta$ and $\S$,  cf. \cite{GZDC}:
$$\frac{\Delta(t)}{1-t}=\sum_{s\in \S}\ t^s.$$
Since $\Delta(1)=1$ and $\Delta'(1)=\delta$, one gets
$$\Delta(t)=1+\delta(t-1)+(t-1)^2\cdot Q(t)$$
for some polynomial $Q(t)=\sum _{i=0}^{2\delta-2} \alpha_it^i$ of
degree $2\delta-2$ with integral coefficients. In fact, all the
coefficients $\{\alpha_i\}_{i}$ are strictly positive,
and:
$$\delta=\alpha_0\geq \alpha_1\geq \cdots \geq \alpha_{2\delta-2}=1.$$
Indeed, by the above identity
$\delta+(t-1)Q(t)=\sum_{s\not\in\S}t^s$, or
$Q(t)=\sum_{s\not\in\Gamma}(t^{s-1}+\cdots +t+1)$. Therefore
$$\alpha_i=\#\{s\not\in\S\ :\ s>i\}.$$
Sometimes it is  convenient
to consider the symmetrized Alexander polynomial $\Delta^\#(t)=t^{-\delta}\Delta(t)$ too. Its second derivative
at 1 satisfies
$$(\Delta^\#)''(1)=2Q(1)+\delta-\delta^2.$$

\vspace{2mm}

Regarding Heegaard Floer homologies  we will use the following notations, cf. \cite{OSzP,Nem}.
Consider the $\Z[U]$--module $\Z[U,U^{-1}]$, and  denote by
 $\calt_0^+$ its quotient by the submodule  $U\cdot \Z[U]$.
It is a $\Z[U]$--module with grading  $\deg(U^{-h})=2h$.
Similarly, for any $n\geq 1$, define the $\Z[U]$--module
$\calt_0(n)$ as the quotient of $\Z\langle U^{-(n-1)},
U^{-(n-2)},\ldots, 1,U,\ldots \rangle$ by $U\cdot \Z[U]$ (with
the same grading). Hence, $\calt_0(n)$, as a $\Z$--module, is the
free $\Z$--module $\Z \langle 1,U^{-1},\ldots,U^{-(n-1)} \rangle$
with  finite  $\Z$--rank $n$.

More generally,  fix an arbitrary $\Q$-graded $\Z[U]$-module $P$
with $h$-homogeneous elements $P_h$. Then for any $x\in \Q$  we
denote by $P[x]$ the same module graded in such a way that
$P[x]_{h+x}=P_{h}$. Then set $\calt^+_x:=\calt^+_0[x]$ and
$\calt_x(n):=\calt_0(n)[x]$.

\vspace{2mm}

Now we are ready to prove
Theorem \ref{thmintr:HF} regarding the integral Heegaard Floer homology  and the
Casson (or Seiberg--Witten) invariant
  of $-\Sigma=S^3_1(T_{p.q})$, cf. Remark \ref{rem:intr}.
\begin{proof}[Proof of Theorem \ref{thmintr:HF}]
The statements are direct consequences of Proposition \ref{prop:minoftau} and Theorem 8.3 of \cite{Nem}.
We only have to check the numerical data; for this see also Proposition 3.5.2 and Corollary 3.7 of \cite{Nem}
for the lattice cohomology associated with graded roots.

In the second identity the torsion modules appear symmetrically, their lengths are respectively
$\tau(M_{\delta-2-k})-\tau(m_{\delta-2-k})$ and $\tau(M_{\delta-1+k})-\tau(m_{\delta+k})$
($0\leq k\leq \delta-2$), both equal to $\alpha_{\delta+k}$. The weight of elements in ${\rm ker}\,U$ are
$2\tau(m_{\delta-2-k})-(K^2+s)/4=2\tau(M_{\delta-2-k})-2\alpha_{\delta+k}+\delta^2-3\delta=k^2+k-2\alpha_{\delta+k}$.
The third identity is exactly this identity for $k=-1$; it is the statement of Theorem \ref{thm:main} as well.
For the fourth identity we use either the second one combined with the identities
$\alpha_i-\alpha_{2\delta-2-i}=\delta-i-1$,  or, directly Corollary 3.7 of \cite{Nem}. 
\end{proof}

The above results --- namely \eqref{eq:K}, Proposition \ref{prop:minoftau} and Theorem \ref{thmintr:HF} --- bear strong
resemblance to the case of negative surgeries  \cite{Nem-neg}. We now cite these results.

\begin{proposition}\label{prop:s'}
 Consider the negative definite Seifert 3--manifold $\Sigma':=S^3_{-1}(T_{p,q})$.
Then one has the following facts:

(a) $K^2+s$ and the $\tau$ function associated with  $\Sigma'$ satisfy:
$$(K^2+s)(\Sigma')=-4\delta(\delta-1).$$
$\tau$  has $2\delta-1$ local maxima at $M_n$ $(0\leq n\leq 2\delta-2$) and
$2\delta$ local minima at $m_n$ ($0\leq n\leq 2\delta-1$), where
\begin{align*}
&\tau(m_n)=\frac{n(n-2\delta+1)}{2}, \ \ \ \min\,\tau=\tau(m_{\delta-1})=\tau(m_\delta)=-\delta(\delta-1)/2,\\
&\tau(M_n)-\tau(m_n)=\alpha_{2\delta-2-n},\\
&\tau(M_n)-\tau(m_{n+1})=\alpha_n,
\end{align*}

(b) The Heegaard Floer homology of $-\Sigma'$ satisfies:
\begin{align*}
&HF^+_{odd}(-\Sigma')=0,\ \ \ \ d(-\Sigma')=0,\\
&HF^+_{even}(-\Sigma')=\calt^+_{0}\oplus \bigoplus_{k=0}^{\delta-2}
\calt_{(k+1)(k+2)}(\alpha_{\delta+k})^{\oplus  2}\oplus \calt_0(\alpha_{\delta-1}),\\
&\lambda(-\Sigma')=\mathbf{sw}(-\Sigma')={\rm rank}\,
HF^+_{red}(-\Sigma')=Q(1)-\frac{\delta(\delta-1)}{2}=\frac{1}{2}(\Delta^{\#})''(1).\end{align*}
\end{proposition}

In the case of $\Sigma'=S^3_{-1}(T_{p,q})$
all the local minima are easily determinable numbers depending only on $\delta$, while  the
local maxima depend essentially on the coefficients $\alpha_i$.
In contrast, for $\Sigma=-S^3_{1}(T_{p,q})$
all the local maxima  depend only on $\delta$, while  the
local minima  on the coefficients $\alpha_i$.

In order to explain more deeply this parallelism, we consider the following construction.
Assume that we have  a sequence of $v+1$ integers $\{\beta_0,\ldots,\beta_v\}$ with
$\beta_0\geq \beta_1\geq \cdots \geq \beta_v\geq 1$. Then, using this sequence one can construct a
continuous sawtooth diagram with $v$ teeth --- as the graph of a function $\tau:[0,2v+2]\to\R$ --- as follows.

We begin with $\tau(0)=0$, then $\tau$ is affine on each interval $[\beta_i,\beta_{i+1}]$; and for each $i$ 
it increases
by $\beta_{v-i}$ on the interval $[2i,2i+1]$, while it decreases by $\beta_i$ on  $[2i+1,2i+2]$.

Let us now consider the coefficients $\{\alpha_0,\alpha_1,\ldots,\alpha_{2\delta-2}\}$ of $Q(t)$ associated with
the semigroup $\S$, or with the Alexander polynomial $\Delta$.
Then the $\tau$ function of $\Sigma'= S^3_{-1}(T_{p,q})$ is associated exactly with this sequence
$\{\alpha_0,\alpha_1,\ldots,\alpha_{2\delta-2}\}$ (and it has in the Heegaard Floer theory
the normalization term $K^2+s=-4\delta(\delta-1)$),
while the $\tau$ function of $\Sigma= -S^3_{1}(T_{p,q})$ is associated  with the shorter sequence
$\{\alpha_1,\ldots,\alpha_{2\delta-2}\}$, hence it has one tooth less (and  $K^2+s=-4\delta(\delta-3)$).
See Figure 1 for $(p,q)=(3,4)$.

\begin{figure}
\begin{pspicture}(-6,-3.5)(6,1)
\psline(-6,0)(-5.5,0.5)(-5,-1)(-4.5,-0.5)(-4,-1.5)(-3.5,-1)(-3,-1.5)(-2.5,-0.5)(-2,-1)(-1.5,0.5)(-1,0)
\psline(1,0)(1.5,0.5)(2,-0.5)(2.5,0)(3,-0.5)(3.5,0)(4,-0.5)(4.5,0.5)(5,0)
\rput(-5.75,-2.3){$1$}
\rput(-5.25,-2.8){$3$}
\rput(-4.75,-2.3){$1$}
\rput(-4.25,-2.8){$2$}
\rput(-3.75,-2.3){$1$}
\rput(-3.25,-2.8){$1$}
\rput(-2.75,-2.3){$2$}
\rput(-2.25,-2.8){$1$}
\rput(-1.75,-2.3){$3$}
\rput(-1.25,-2.8){$1$}
\rput(1.25,-2.3){$1$}
\rput(1.75,-2.8){$2$}
\rput(2.25,-2.3){$1$}
\rput(2.75,-2.8){$1$}
\rput(3.25,-2.3){$1$}
\rput(3.75,-2.8){$1$}
\rput(4.25,-2.3){$2$}
\rput(4.75,-2.8){$1$}
\psline[linestyle=dotted](-6.2,0)(-0.8,0)
\psline[linestyle=dotted](0.8,0)(5.5,0)
\rput(6.15,0){$0$}
\psline[linestyle=dotted](-6.2,0.5)(-0.8,0.5)
\psline[linestyle=dotted](0.8,0.5)(5.5,0.5)
\rput(6.15,0.5){$1$}
\psline[linestyle=dotted](-6.2,-0.5)(-0.8,-0.5)
\psline[linestyle=dotted](0.8,-0.5)(5.5,-0.5)
\rput(6.0,-0.5){$-1$}
\psline[linestyle=dotted](-6.2,-1)(-0.8,-1)
\psline[linestyle=dotted](0.8,-1)(5.5,-1)
\rput(6.0,-1){$-2$}
\psline[linestyle=dotted](-6.2,-1.5)(-0.8,-1.5)
\psline[linestyle=dotted](0.8,-1.5)(5.5,-1.5)
\rput(6.0,-1.5){$-3$}

\psline[linestyle=dashed](-6,0)(-6,-3)
\psline[linestyle=dashed](-5.5,0.5)(-5.5,-3)
\psline[linestyle=dashed](-5,-1)(-5,-3)
\psline[linestyle=dashed](-4.5,-0.5)(-4.5,-3)
\psline[linestyle=dashed](-4,-1.5)(-4,-3)
\psline[linestyle=dashed](-3.5,-1)(-3.5,-3)
\psline[linestyle=dashed](-3,-1.5)(-3,-3)
\psline[linestyle=dashed](-2.5,-0.5)(-2.5,-3)
\psline[linestyle=dashed](-2,-1)(-2,-3)
\psline[linestyle=dashed](-1.5,0.5)(-1.5,-3)
\psline[linestyle=dashed](-1,0)(-1,-3)

\psline[linestyle=dashed](1.0,0.0)(1.0,-3)
\psline[linestyle=dashed](1.5,0.5)(1.5,-3)
\psline[linestyle=dashed](2.0,-0.5)(2.0,-3)
\psline[linestyle=dashed](2.5,0.0)(2.5,-3)
\psline[linestyle=dashed](3.0,-0.5)(3.0,-3)
\psline[linestyle=dashed](3.5,0.0)(3.5,-3)
\psline[linestyle=dashed](4.0,-0.5)(4.0,-3)
\psline[linestyle=dashed](4.5,0.5)(4.5,-3)
\psline[linestyle=dashed](5.0,0.0)(5.0,-3)

\pscircle[fillcolor=black,fillstyle=solid](-6,0){0.06}
\pscircle[fillcolor=black,fillstyle=solid](-5.5,0.5){0.06}
\pscircle[fillcolor=black,fillstyle=solid](-5,-1){0.06}
\pscircle[fillcolor=black,fillstyle=solid](-4.5,-0.5){0.06}
\pscircle[fillcolor=black,fillstyle=solid](-4,-1.5){0.06}
\pscircle[fillcolor=black,fillstyle=solid](-3.5,-1){0.06}
\pscircle[fillcolor=black,fillstyle=solid](-3,-1.5){0.06}
\pscircle[fillcolor=black,fillstyle=solid](-2.5,-0.5){0.06}
\pscircle[fillcolor=black,fillstyle=solid](-2,-1){0.06}
\pscircle[fillcolor=black,fillstyle=solid](-1.5,0.5){0.06}
\pscircle[fillcolor=black,fillstyle=solid](-1,0){0.06}
\pscircle[fillcolor=black,fillstyle=solid](1,0){0.06}
\pscircle[fillcolor=black,fillstyle=solid](1.5,0.5){0.06}
\pscircle[fillcolor=black,fillstyle=solid](2,-0.5){0.06}
\pscircle[fillcolor=black,fillstyle=solid](2.5,0){0.06}
\pscircle[fillcolor=black,fillstyle=solid](3,-0.5){0.06}
\pscircle[fillcolor=black,fillstyle=solid](3.5,0){0.06}
\pscircle[fillcolor=black,fillstyle=solid](4,-0.5){0.06}
\pscircle[fillcolor=black,fillstyle=solid](4.5,0.5){0.06}
\pscircle[fillcolor=black,fillstyle=solid](5,0){0.06}

\end{pspicture}
\caption{The $\tau$ function for $S^3_{-1}(T_{3,4})$ on the left and $-S^3_1(T_{3,4})$ on the right. Below we show the increment and the decrement
of the $\tau$ function on each interval.  $S_{3,4}=\mathbb{N}\setminus\{1,2,5\}$,
hence  $(\alpha_0,\alpha_1,\alpha_2,\alpha_3,\alpha_4)=(3,2,1,1,1)$.}
\end{figure}

From these data all the numerical properties of the corresponding Heegaard Floer homologies follow
by the general theory of graded roots, see \cite{Nem}.
Note that in the case of $\Sigma$, the missing entry is $\alpha_0=\delta$, hence once we know $\delta$,
the two sequences   determine each other!

We expect that this `duality' is valid in a  more general situation.

\section{Relation to the Levine--Tristram signature.}\label{section:TL}

Let $p$ and $q$ be  two coprime integers as above. For $y\in[0,1]$  we denote by $\sigma(y)$
the Levine--Tristram signature of the (positive) torus knot $T_{p,q}$ evaluated at $e^{2\pi iy}$.

Moreover, let $Sp:=\{\frac{i}{p}+\frac{j}{q}\colon 1\le i\le p-1;1\le j\le q-1\}$ be the spectrum of the
local plane curve singularity $x^p+y^q=0,\ (x,y)\in ({\mathbb C}^2,0)$.

\begin{proposition}\label{prop:h-int} With the above notations,
for any integer $a\in [0,pq)$  one has
 \[\#\{s\not\in\S\,:\, s\geq a\}=\#\,\{\,[1+\frac{a}{pq},2)\,\cap\,Sp\},
 \]
and
\begin{equation}\label{eq:h-int}
\#\{s\not\in S_{p,q}\,:\, s\ge a\}=\delta+\frac14\cdot \sigma\big(\frac{a}{pq}\big)
-\frac12\big(a-\intfrac{a}{p}-\intfrac{a}{q}-\tilde{c}(a)\big),
\end{equation}
where
\[
\tilde{c}(a)=\left\{\begin{array}{ll}
\frac12 & \ \mbox{if \ $1+\frac{a}{pq}\in Sp$},\\
-\frac12 & \ \mbox{if \ $\frac{a}{pq}\in Sp$},\\
0 & \ \mbox{otherwise.}\end{array}\right.
\]
\end{proposition}
\begin{proof} If $s\not\in \S$ and $s\geq a$, then $s+pq=\alpha p+\beta q$ for some $\alpha\in[0,q)$ and
$\beta \in[0,p)$ (see e.g. \eqref{eq:sgr}).
But $\alpha\beta\not=0$ (otherwise $s\in \S$), hence $1+\frac{s}{pq}\in Sp\cap [1+\frac{a}{pq},2)$.
This proves the first identity.

Let us define the following four quantities $S_i:=\#\,Sp\cap I_i$ ($1\leq i\leq 4$), where
\begin{equation*}
I_1=(0,\frac{a}{pq}], \ \ I_2=(\frac{a}{pq},1), \ \
I_3=(1,1+\frac{a}{pq}), \ \  I_{4}=[1+\frac{a}{pq},2).
\end{equation*}
Then  $S_1+S_2=S_3+S_4=\delta$. Moreover, by \cite[Corollary~4.4.9]{BN}  (see also  \cite{Lith}),
one has
\[\sigma\big(\frac{a}{pq}\big)=S_1+S_4-S_2-S_3+\#(Sp\cap\{\frac{a}{pq},1+\frac{a}{pq}\}).\]
On the other hand, we claim that
\[S_1+S_3+ \#(Sp\cap\{1+\frac{a}{pq}\})=a-\intfrac{a}{p}-\intfrac{a}{q}.\]
Indeed, the left hand side is $\big(\,(0,\frac{a}{pq}]\cup (1,1+\frac{a}{pq}]\,\big)\cap Sp$. Then
$\frac{s}{pq}$ is in this set if and only if
$s$ is an integer in $[0,x]$ which is not divisible either by  $p$ or by $q$. Their number is the right hand side
of the identity.

In order to end the proof, we write
\[4S_4=S_1+S_4-S_3-S_2+3(S_4+S_3)+(S_2+S_1)-2(S_1+S_3)\]
and substitute the above identities.
\end{proof}

Substituting $a=\delta$ into \eqref{eq:h-int} we obtain
\begin{corollary}\label{cor:sigma}
$$d(S^3_1(T_{p,q}))=-2
\#\{s\not\in S_{p,q}\,:\, s\ge \delta\}=-\delta-\frac12\cdot \sigma\big(\frac{\delta}{pq}\big)
-\intfrac{\delta}{p}-\intfrac{\delta}{q}-\tilde{c}(\delta).$$
\end{corollary}


\section{Inequalities}\label{section:INEQ}

Fix  $p<q$.
Consider the surface Brieskorn singularity 
$(u,v,w)\in(\mathbb{C}^3,0)\colon u^p+v^q+w^2=0$, i.e. the double suspension of the plane curve
singularity $x^p+y^q=0$.
Its Milnor number is $\mu=2\delta$; let $\mu_{+}$ and $\mu_{-}$ be the dimensions of maximal subspaces of the
vanishing homology where the intersection form is positive/negative definite. (Note that the intersection form
is non--degenerate, hence the dimension $\mu_0$ of its kernel is zero.) Therefore, $2\delta=\mu_++\mu_-$
and the signature is $\sigma=\sigma(K)=\mu_+-\mu_-$.

\begin{lemma}\label{lemma:INEQ} For $p$ and $q$ relative primes one has:

\vspace{1mm}

(a) \ $\#\{s\not\in\S\,:\, s\geq \delta\}\,\leq \, \frac{\delta+1}{2}$

\vspace{1mm}

(b) \ $\delta+1\leq -\sigma$.

\vspace{1mm}

The inequality (a) is sharp for $p=2$ and $\delta$ odd, while (b) is equality for $(p,q)=(2,3)$.
\end{lemma}
\begin{proof}
(a) Set $S^*:=\{s\in\S\,:\, 0<s<\delta\}$, and let $s^*$ be its maximal element. Then
$\{0\}\sqcup S^*\sqcup (s^*+S^*)\subset \S\cap [0,2\delta-1)$, hence
$1+2\#S^*\leq \#[0,2\delta-1)-\#({\mathbb N}\setminus \S)=\delta$. Therefore,
$\#\{s\not\in\S\,:\, s\geq \delta\}=\#\{s\in\S\,:\, s\leq \delta-1\}=1+\#S^*\leq (\delta+1)/2$.

(b) By Theorem 4.1 of \cite{N}, $-6\sigma\geq 6\delta+q-1$, hence the inequality follows for
$q\geq 7$. For the other pairs when $p<q\leq 6$ can be checked case by case.
\end{proof}
Since the four genus $g_4$ of $T_{p,q}$ is $\delta$, we get the statement  of Corollary \ref{cor:INEQ}(a).

Having the equality $-d(K)=2\ceilfrac{-\sigma(K)}{4}$ for alternating knots, we might wonder
if for any torus knot $-d\leq 2\ceilfrac{-\sigma}{4}$ is valid.  However, this is not the case:
for the pair $(p,q)=(4,5)$ one has $-d=6$ and $\sigma=-8$.

Nevertheless, asymptotically for `most' of the pairs $(p,q)$, $-d$ grows like $-\sigma/2$.
\begin{lemma}
Recall that $p<q$. Then
$$\#Sp\cap \big[1+\frac{\delta}{pq}\big)\leq q-1+\#Sp\cap\big[\frac32,2\big).$$
\end{lemma}
\begin{proof}
We need to show that $\#Sp\cap [1+\delta/pq,3/2)\leq q-1$. For this notice that if
$s=\alpha/p+\beta/q$ is in the interval $[1+\delta/pq,3/2)$, then $s+\frac{1}{p}\geq 3/2$.
Hence, $1\leq \beta\leq q-1$ determines any $s$ in that interval.
\end{proof}
Note that $\#Sp\cap [3/2,2)=\mu_+/2$ by Thom--Sebastiani type theorem for the spectrum of the
suspension and by relationship between the spectral number of surfaces and the intersection form \cite{St}.
Note also that the inequality $\delta\leq -\sigma-1$ of \ref{lemma:INEQ}(b) can be rewritten as
$\mu_+\leq (\delta-1)/2$. Therefore, we obtain:
\begin{equation}\label{eq:last}
-d(S^3_1(T_{p,q}))\leq 2q-2+\mu_+\leq 2q-2+\frac{\delta-1}{2}.
\end{equation}
If $p$ is `small' with respect to $q$, then the $q$ term  at the right hand side makes the inequality
weak. Nevertheless, if $p$ `grows together with $q$', then one can find a positive constant
$c$ such that $2q-\frac52\leq c\sqrt{\delta}$.
For example, if $p=q-1$ then $2q-\frac52\leq 8\sqrt{\delta}+1$.
Hence \eqref{eq:last} together with Lemma \ref{lemma:INEQ}(b) provides Corollary \ref{cor:INEQ}(b).

\begin{acknowledgements}
The first author is supported by Polish MNiSzW Grant No N N201 397937 and also by a Foundation for Polish
Science FNP (program KOLUMB). He also wishes to express his thanks to Renyi Insitute of Mathematics for hospitality.
The second author is partially supported by OTKA Grant K67928.
\end{acknowledgements}

\end{document}